\documentclass[12pt]{article}
\usepackage[centertags]{amsmath}
\usepackage{amsfonts}
\usepackage{amssymb}
\usepackage{amsthm}
\usepackage{newlfont}
\usepackage[left=1in,top=1in,right=1in,nohead,foot=0.5in]{geometry}
 \newtheorem{thm}{Theorem}[section]
 \newtheorem{cor}[thm]{Corollary}
 \newtheorem{lem}[thm]{Lemma}
 \newtheorem{prop}[thm]{Proposition}
 \theoremstyle{definition}
 
 \theoremstyle{remark}
 \newtheorem{rem}[thm]{Remark}
 \numberwithin{equation}{section}
 \newcommand{\Spec}{\operatorname{Spec}}
 \newcommand{\Proj}{\operatorname{Proj}}
 \newcommand{\Res}{\operatorname{Res}}
 \newcommand{\Hom}{\operatorname{Hom}}
 \newcommand{\Rat}{\operatorname{Rat}}
 \newcommand{\GL}{\operatorname{GL}}
 \newcommand{\SL}{\operatorname{SL}}
 \newcommand{\PGL}{\operatorname{PGL}}
 \newcommand{\Stab}{\operatorname{Stab}}
 
 \newcommand{\Fix}{\operatorname{Fix}}

 \newcommand{\chr}{\operatorname{char}}
 
\newlength{\defbaselineskip}
\begin{document}

\title{The Space of Morphisms on Projective Space}
\author{Alon Levy}

\maketitle

\abstract{The theory of moduli of morphisms on $\mathbb{P}^{n}$ generalizes the study of rational maps on $\mathbb{P}^{1}$. This paper proves three results about
the space of morphisms on $\mathbb{P}^{n}$ of degree $d > 1$, and its quotient by the conjugation action of $\PGL(n+1)$. First, we prove that this quotient is
geometric, and compute the stable and semistable completions of the space of morphisms. This strengthens previous results of Silverman, as well as of Petsche,
Szpiro, and Tepper. Second, we bound the size of the stabilizer group in $\PGL(n+1)$ of every morphism in terms of only $n$ and $d$. Third, we specialize to the
case where $n = 1$, and show that the quotient space is rational for all $d > 1$; this partly generalizes a result of Silverman about the case $d = 2$.

\section{Introduction and Notation}

\noindent A rational map from $\mathbb{P}^{n}$ to itself is determined by an $(n+1)$-tuple of polynomials in $n+1$ variables, all homogeneous of the same degree
$d$. If this map is a morphism, it will be finite of degree $d^{n}$. In the rest of this paper, we will refer to such a rational map as a degree $d$ map on
$\mathbb{P}^{n}$ by abuse of notation. The space of degree $d$ maps on $\mathbb{P}^{n}$ is projective, with homogeneous coordinates coming from monomials of degree
$d$. There are ${n+d \choose d}$ such monomials, so that this space has dimension ${n+d \choose d}(n+1) - 1$. We write $N_{d}^{n}$ for the dimension of this space,
or $N$ when $d$ and $n$ are clear.

\medskip The case of interest is morphisms on $\mathbb{P}^{n}$. In the sequel, we refer to the polynomials defining the map as $q_{0}, q_{1}, \ldots, q_{n}$. Then a
map $(q_{0}:\ldots:q_{n})$ is a morphism if and only if the $q_{i}$'s share no common geometric root. The $q_{i}$'s only share a common root on a hypersurface of
$\mathbb{P}^{N}$ which we call the resultant subvariety and which is defined over $\mathbb{Z}$; we denote its complement by $\Hom_{d}^{n}$.

\medskip The space $\mathbb{P}^{N}$ of rational maps comes equipped with an action of $\PGL(n+1)$ by conjugation. The conjugation action $A\cdot\varphi = A\varphi
A^{-1}$, fixes the resultant, which gives an action of $\PGL(n+1)$ on $\Hom_{d}^{n}$. In this paper, we mainly study the quotient of this action, which we denote
$\mathrm{M}_{d}^{n}$, or $\mathrm{M}_{d}$ when $n = 1$. We will show that this quotient is geometric in the sense of geometric invariant theory \cite{GIT}, and
compute the largest stable and semistable loci $\Hom_{d}^{n, s}$ and $\Hom_{d}^{n, ss}$, which satisfy $\Hom_{d}^{n} \subset \Hom_{d}^{n, s} \subset \Hom_{d}^{n,
ss} \subset \mathbb{P}^{N}$.

\medskip Knowing that the quotient $\mathrm{M}_{d}^{n}$ is well-behaved is often necessary to answer questions about the geometry of families of dynamical systems.
In \cite{PST}, Petsche, Szpiro, and Tepper prove that $\mathrm{M}_{d}^{n}$ exists as a geometric quotient in order to show that isotriviality is equivalent to
potential good reduction for morphisms of $\mathbb{P}^{n}$ over function fields, generalizing previous results in the one-dimensional case. In \cite{DeM07}, DeMarco
uses the explicit description of the space $\mathrm{M}_{2}$ in order to study iterations of quadratic maps on $\mathbb{P}^{1}$, and one can expect similar results
in higher dimension given a better understanding of the structure of $\mathrm{M}_{d}^{n}$.

\medskip By now the theory of morphisms on $\mathbb{P}^{1}$ is the standard example in dynamical systems. For a survey of the arithmetic theory, see \cite{ADS};
also see a recent paper by Manes \cite{Man} about moduli of morphisms on $\mathbb{P}^{1}$ with a marked point of period $n$, which functions as a dynamical level
structure. In the complex case, see an overview by Milnor \cite{D1C}, and the work of DeMarco \cite{DeM05} \cite{DeM07} about compactifications of the space
$\mathrm{M}_{d}$ that respect the iteration map. Despite this, the higher-dimensional theory remains understudied. The only prior result in the direction of moduli
of morphisms on $\mathbb{P}^{n}$ is the proof in \cite{PST} that $\mathrm{M}_{d}^{n}$ exists as a geometric quotient. Unfortunately, the proof does not lend itself
well to finding the stable and semistable spaces for the action of $\PGL(n+1)$ on $\mathbb{P}^{N}$, nor does it bound the size of the finite stabilizer group
uniformly on $\Hom_{d}^{n}$.

\medskip The first two tasks in this paper are then to construct alternative proofs of the fact that the quotient $\mathrm{M}_{d}^{n}$ is geometric, first by
explicitly describing the stable and semistable loci, and second by finding a uniform bound for the size of the stabilizer group in $\PGL(n+1)$. The former we will
do in section $2$, using the Hilbert-Mumford criterion for stability and semistability. We will see that the complements of both $\Hom_{d}^{n, s}$ and $\Hom_{d}^{n,
ss}$ are equal to a finite union of linear subvarieties and their $\PGL(n+1)$-conjugates; this contrasts with the $n = 1$ case, when the complement is the
$\PGL(2)$-orbit of only one linear subvariety. In section $3$ we will study the stabilizer groups, proving a uniform bound on their sizes, valid over all fields and
rings of definition, depending only on $n$ and $d$. This will strengthen previous results in this direction for $n = 1$ in \cite{Sil94}.

\medskip Most results in this paper are a natural generalization of the study of morphisms on $\mathbb{P}^{1}$ in \cite{Sil96}, which refers to the space of
morphisms as $\Rat_{d}$ and its quotient as $\mathrm{M}_{d}$, and which proves that $\mathrm{M}_{2} \cong_{\Spec\mathbb{Z}} \mathbb{A}^{2}$ using the theories of
fixed points and multipliers. Specializing to the case where $n = 1$, we will prove that $\mathrm{M}_{d}$ is rational for all $d$ in section $4$. This is new even
in the case of $d = 3$. The proof in this paper is based on showing that $\mathrm{M}_{d}$ is birational to a vector bundle over the space $\mathrm{M}_{0, d+1}$ of
$d + 1$ unmarked points on $\mathbb{P}^{1}$, which is known to be rational.

\medskip Unfortunately, we do not see any easy generalization of rationality to $\mathrm{M}_{d}^{n}$. The obstruction is that the space of unmarked points on
$\mathbb{P}^{n}$ is not known to be rational. Clearly $\Hom_{d}^{n}$ is rational, so $\mathrm{M}_{d}^{n}$ is unirational, which for some applications, such as the
density of points defined over a number field $K$, is enough. However, in order to investigate the structure of $\mathrm{M}_{d}^{n}$ we need more than that. We do
not expect a result along the lines of that in \cite{Sil96}, that $\mathrm{M}_{2} \cong \mathbb{A}^{2}$, but we do expect rationality of $\mathrm{M}_{d}^{n}$.

\medskip I would like to express my gratitude to my advisor Shouwu Zhang for introducing me to dynamical systems and guiding my research, to Lucien Szpiro and Joe
Silverman for looking at the proofs of the major theorems in this paper, and to Xander Faber for helping me with this paper's presentation.

\section{The Spaces $\Hom_{d}^{n}$ and $\mathrm{M}_{d}^{n}$}

\noindent The space $\Hom_{d}^{n}$ of degree-$d$ morphisms on $\mathbb{P}^{n}$ arises as the subset of $\mathbb{P}^{N} = \{(q_{0}:q_{1}:\ldots:q_{n})\}$ defined by
the condition that the $q_{i}$'s share no common root. In order to give this space an algebraic structure, we investigate its complement. We will show the following
result, proven by Macaulay \cite{Mac} and reinterpreted here in modern language:

\begin{thm}\label{Macaulay}The maps on $\mathbb{P}^{n}$ of degree $d$ such that the $q_{i}$'s share a nonzero root form a closed, irreducible subvariety of
$\mathbb{P}^{N}$ of codimension $1$, which is defined over $\mathbb{Z}$.\end{thm}

\begin{proof}Consider the variety $V = \mathbb{P}^{n} \times \mathbb{P}^{N}$. We think of $V$ as representing a set of polynomials $(q_{0}: q_{1}: \ldots: q_{n})$
acting on the point $(x_{0}: x_{1}: \ldots: x_{n})$. Consider the resultant subvariety $U \subset V$ defined by the condition that $q_{i}(\mathbf{x}) = 0$ for all
$i$. This variety clearly has codimension at most $n+1$. If we denote the variables defining $\mathbb{P}^{N}$ as $a^{i}_{j^{i}_{0}j^{i}_{1}\ldots j^{i}_{n}}$ with
$j^{i}_{0} + \ldots + j^{i}_{n} = d$, representing the $x_{0}^{j^{i}_{0}}\ldots x_{n}^{j^{i}_{n}}$ monomial of $q_{i}$, then we see that $U$ is defined by equations
that are bihomogeneous of degree $1$ in the $a^{i}_{J}$'s and $d$ in the $x_{i}$'s.

\medskip We claim that $U$ is irreducible. The claim follows from a generalization of the fact that a primitive polynomial is irreducible over a domain whenever it
is irreducible over its fraction field. More precisely, let $R$ be a domain with fraction field $K$, and let $I$ be an ideal of $R[y_{1}, \ldots, y_{m}]$ that is
not contained in any prime of $R$. We have a natural map $f$ from $\Spec K[y_{1}, \ldots, y_{m}]$ to $\Spec R[y_{1}, \ldots, y_{m}]$. If $V(I)$ is reducible over
$R$, say $V(I) = V_{1} \cup V_{2}$ with $V_{i}$ nonempty, then either $V(I)$ is reducible over $K$, or one $f^{-1}(V_{i})$, say $f^{-1}(V_{1})$, is empty. In the
latter case, $I(V_{1})$ may not contain nonconstant polynomials, so it contains at least one prime constant. This contradicts the assumption that $I$ is not
contained in any prime of $R$; hence, $V(I)$ is reducible over $K$.

\medskip With the above generalization, suppose that $U$ is reducible. Then it is also reducible as a subvariety of $\mathbb{A}^{n+1}\times\mathbb{A}^{N+1}$.
Further, by letting $R = \mathbb{Z}[x_{0}, \ldots, x_{n}]$ and $K$ be its fraction field, we see that either $U$ is contained in a prime of $R$, or $U$ is reducible
in $\mathbb{A}^{N+1}_{K}$. The former case is impossible since $U$ is not contained in any prime of $\mathbb{Z}$ or any relevant prime ideal of the ring of
polynomials over $\mathbb{Z}$, and the latter is impossible since it is defined by linear equations in the $a^{i}_{J}$'s. Either way this is a contradiction, so $U$
is irreducible and the claim is proven.

\medskip Finally, the maps on $\mathbb{P}^{n}$ of degree $d$ whose polynomials have a common nonzero root arise as the projection of $U$ onto the second factor of
$\mathbb{P}^{n} \times \mathbb{P}^{N}$. It is irreducible because the projection map is surjective. It is closed because the map is proper. It has codimension at
most $1$ because almost all polynomials in $U$ share just one root, so that the dimension of $U$ and its image are equal. It has exact codimension $1$ because some
maps, for instance $q_{i} = x_{i}^{d}$, are morphisms. And it is defined over $\mathbb{Z}$ because every construction we have made in this proof is defined over
$\mathbb{Z}$.\end{proof}

\medskip We call the image of $U$ the resultant subvariety of $\mathbb{P}^{N}$; we call its generating polynomial the Macaulay resultant and denote it by
$\Res_{d}^{n}$. Macaulay proved the theorem by constructing the resultant explicitly, and showing that it has integer coefficients and is irreducible. His explicit
construction shows that if the polynomials are homogeneous of degrees $d_{0}, d_{1}, \ldots, d_{n}$, then the resultant is $(n+1)$-homogeneous in the coefficients
of each polynomial $p_{i}$ of degree $\prod_{j \neq i}d_{j}$. In our case, all the degrees are equal to $d$, so that the resultant is $(n+1)$-homogeneous in the
coefficients of each $q_{i}$ of degree $d^{n}$. In particular, the resultant subvariety is a hypersurface of degree $(n+1)d^{n}$.

\medskip Theorem~\ref{Macaulay} shows that the space of morphisms is the complement of the resultant subvariety, and is therefore affine and of dimension $N$.
Silverman \cite{Sil96}, who only considers the case $n = 1$, refers to this space as $\Rat_{d}$; we will refer to it as $\Hom_{d}^{n}$ and to its complement in
$\mathbb{P}^{N}$ as $\Res_{d}^{n}$ by abuse of notation.

\medskip The action of $\PGL(n+1)$ on $\mathbb{P}^{n}$ leads to a conjugation action on $\Hom_{d}^{n}$, wherein $A \in \PGL(n+1)$ acts on a rational map $\varphi$
by sending it to $A\varphi A^{-1}$. The property of being ill-defined at a point is stable under both the left action mapping $\varphi$ to $A\varphi$ and the right
action mapping $\varphi$ to $\varphi A^{-1}$; hence, the conjugation action is well-defined on $\Hom_{d}^{n}$. The space of endomorphisms of $\mathbb{P}^{n}$
defined by degree-$d$ polynomials may be regarded as the quotient of $\Hom_{d}^{n}$ by the conjugation action.

\medskip \emph{A priori}, we only know that over an algebraically closed field, the quotient exists as a set. In order to give it algebraic structure, we need to
pass to the stable or semistable space in geometric invariant theory \cite{GIT}. Fortunately, we have the following result:

\begin{thm}\label{stable}Every $\varphi \in \Hom_{d}^{n}$ is stable.\end{thm}

\begin{proof}We use the Hilbert-Mumford criterion, as described in chapter $2$ of \cite{GIT}. To do that, we pull back the action of $\PGL(n+1)$ on $\mathbb{P}^{N}$
to the action of $\SL(n+1)$ on $\mathbb{A}^{N+1}$, and consider one-parameter subgroups of $\SL(n+1)$. The criterion states that a point lies in the stable space
$\Hom_{d}^{n, s}$ (respectively, the semistable space $\Hom_{d}^{n, ss}$) iff for every such subgroup, its action on the point can be diagonalized with diagonal
elements $t^{a_{I}}$, and at least one $a_{I}$ is negative (resp. non-positive).

\medskip Note that the action of $A \in \SL(n+1)$ on $\varphi \in \mathbb{A}^{N+1}$ is conjugate to the action of $BAB^{-1}$ on $B\varphi B^{-1}$. In particular, it
will have the same eigenvalues, so the action of a one-parameter subgroup $G = \mathbb{G}_{m}$ will have the same $a_{I}$'s. Therefore, we may conjugate $G$ to be
diagonal, which will be enough to give us criteria for stability and semistability up to conjugation. So from now on, we assume $G$ is the diagonal subgroup whose
$i$th diagonal entry is $t^{a_{i}}, a_{i} \in \mathbb{Z}$. Here we label the rows and columns from $0$ to $n$, in parallel with the label for the $q_{i}$'s. We have
$a_{0} + \ldots + a_{n} = 0$. We may also assume that $a_{0} \geq a_{1} \geq \ldots \geq a_{n}$, after conjugation if necessary, and that the $a_{i}$'s are coprime.

\medskip The action of $G$ on $\mathbb{A}^{N+1}$ is already diagonal. We denote the $\mathbf{x^{d}}$ coefficient of $q_{i}$ by $c_{\mathbf{d}}(i)$; then $G$
multiplies $c_{\mathbf{d}}(i)$ by $t^{a_{i}}t^{-(a_{0}d_{0} + \ldots + a_{n}d_{n})}$. A point $\varphi$ is not stable (resp. unstable) if for some choice of $G$,
all the $c_{\mathbf{d}}(i)$'s for which $a_{0}d_{0} + \ldots + a_{n}d_{n} > a_{i}$ (resp. $a_{0}d_{0} + \ldots + a_{n}d_{n} \geq a_{i}$) are zero. Let us observe
that this means that, for $d > 1$, every $x_{0}^{d}$ coefficient has to be zero, as we will have $da_{0} > a_{0} \geq a_{i}$ for every $i$. This means that
$\varphi$ lacks any $x_{0}^{d}$ coefficient, so that the $q_{i}$'s have a nontrivial zero at $(1:0:\ldots:0)$, and $\varphi \notin \Hom_{d}^{n}$. The property of
not being a morphism is preserved under conjugation, proving the theorem.\end{proof}

\medskip Since $\Hom_{d}^{n}$ is stable, it has a natural geometric quotient induced by the $\PGL(n+1)$ action on $\mathbb{P}^{N}$, which we denote
by $\mathrm{M}_{d}^{n}$; as $\Hom_{d}^{n}$ is affine, $\mathrm{M}_{d}^{n}$ is affine, with structure sheaf $\mathcal{O}_{\Hom_{d}^{n}}^{\SL(n+1)}$. We may also
write $\mathrm{M}_{d}^{n, s}$ for the quotient of the stable space and $\mathrm{M}_{d}^{n, ss}$ for the quotient of the semistable space. The latter quotient is
only categorical, rather than geometric, but will be proper over $\Spec\mathbb{Z}$ (all spaces in question, as well as $\SL(n+1)$, are defined over $\mathbb{Z}$;
hence, so are the quotients).

\medskip Let us now describe the not-stable and unstable spaces more explicitly. In the $n = 1$ case, $G$ depends only on $a_{0}$, which may be taken to be $1$.
This gives us only one criterion for stability (resp. semi-stability), which means that the not-stable (resp. unstable) space is irreducible (in fact, it will be a
linear subvariety and its orbit under $\PGL(2)$-conjugation). When $n > 1$, this is no longer true: $G$ depends on multiple variables, and we can find many infinite
families of coprime $a_{i}$'s that sum to $0$ and are in decreasing order.

\medskip However, the not-stable (resp. unstable) space will still be a union of finitely many linear subvarieties and their $\PGL(n+1)$ conjugates, whose number
will generally grow with $d$ and $n$. This is because there are only $2^{N+1}$ linear spaces defined by conditions of the form $c_{\mathbf{d}}(i) = 0$ for a
collection $J$ of $(\mathbf{d}, i)$ pairs. For each such space, either there exists a $G$ such that $(\mathbf{d}, i) \in J$ if and only if $a_{0}d_{0} + \ldots +
a_{n}d_{n} > a_{i}$ (resp. $a_{0}d_{0} + \ldots + a_{n}d_{n} \geq a_{i}$), or there doesn't. Of course, a given $J$ may correspond to infinitely many $G$, which
will in general have ratios $a_{0}:\ldots:a_{n}$ that are close in the archimedean metric.

\medskip We omit the calculation of the linear subvarieties that occur as the not-stable (resp. unstable) space for each $d$ and $n$, as well as the number of such
varieties. We will just note that there are far fewer than $2^{N+1}$ such varieties: for a start, we have already seen that $((d, 0, \ldots, 0), i) \in J$ for all
$i$. One more constraint that follows trivially from the definition of the $a_{i}$'s is that if $(\mathbf{d}, i) \in J$, then so is $(\mathbf{d}, j)$ for $j > i$.
Put another way, not being stable (resp. instability) imposes more conditions on $q_{j}$ than on $q_{i}$ for $j > i$. It may also be shown that for each $G$ the
number of conditions is roughly between one half and $e^{-1}$ times $N$; we omit the proof, as this result will not be relevant in the remainder of this paper.

\medskip Finally, when $n = 1$, the only $G$ has $a_{0} = 1, a_{1} = -1$, so $a_{0}d_{0} + a_{1}d_{1} = d_{0} - d_{1} = 2d_{0} - d$. When $d$ is even, $2d_{0} - d$
is always even, so the conditions $a_{0}d_{0} + a_{1}d_{1} > a_{i}$ and $a_{0}d_{0} + a_{1}d_{1} \geq a_{i}$ coincide, and the stable and semistable spaces are the
same; this was shown in \cite{Sil96}. We will show that this will never be the case for higher $n$. First, observe that if we set $a_{0} = 1$, $a_{n} = -1$, and
$a_{i} = 0$ for $i \neq 0, n$, we obtain $a_{0}d_{0} + \ldots + a_{n}d_{n} = d_{0} - d_{n}$, which may take any value between $-d$ and $d$ inclusive. Hence, the
conditions $a_{0}d_{0} + \ldots + a_{n}d_{n} > a_{i}$ and $a_{0}d_{0} + \ldots + a_{n}d_{n} \geq a_{i}$ will not coincide.

\medskip Now, suppose that $\varphi$ is a point that is not stable, with $c_{\mathbf{d}}(i) = 0$ if and only if $d_{0} - d_{n} > a_{i}$ with $a_{i}$ as above. If
$\varphi$ is unstable, then we can find some $G$ such that if $a_{0}d_{0} + \ldots + a_{n}d_{n} \geq a_{0}$ then $d_{0} - d_{n} > 1$, and if $a_{0}d_{0} + \ldots +
a_{n}d_{n} \geq a_{i}$ for $i \neq 0, n$, then $d_{0} - d_{n} > 0$. If for that $G$ we have $a_{1} \geq 0$, then looking at the $x_{0}x_{1}^{d-1}$ monomial, we get
$a_{0}d_{0} + \ldots + a_{n}d_{n} = a_{0} + (d-1)a_{1} \geq a_{0}$ but $d_{0} - d_{n} = 1$, a contradiction. If $a_{1} < 0$, then we must have $a_{i} < 0$ for all
$i > 0$, so $a_{0} + a_{n} > 0$. For $d = 2k + 1$, we consider the $x_{0}^{k+1}x_{n}^{k}$ monomial, for which $a_{0}d_{0} + \ldots + a_{n}d_{n} = k(a_{0} + a_{n}) +
a_{0} > a_{0}$ but $d_{0} - d_{n} = 1$; for $d = 2k$, we consider the $x_{0}^{k}x_{n}^{k}$ monomial, for which $a_{0}d_{0} + \ldots + a_{n}d_{n} = k(a_{0} + a_{n})
> 0 > a_{1}$ but $d_{0} - d_{n} = 0$. Either way, we have a contradiction, so $\varphi$ is semistable but not stable. This proves:

\begin{prop}\label{semistable}For all $d, n > 1$, we have $\Hom_{d}^{n, s} \subsetneq \Hom_{d}^{n, ss}$.\end{prop}

\medskip We will conclude this section with the following strict containment:

\begin{prop}\label{containment}$\Hom_{d}^{n} \subsetneq \Hom_{d}^{n, s}$.\end{prop}

\begin{proof}Observe that the linear subvarieties defined above are invariant under conjugation by every upper triangular matrix, at least when we ensure $a_{0}
\geq a_{1} \geq \ldots \geq a_{n}$. Hence, the codimension of the not-stable space is equal to the codimension of the largest linear subvariety, minus
$\frac{n(n+1)}{2}$. It suffices to show this codimension is more than $1$, or, in other words, that every linear subvariety has codimension at least
$\frac{n(n+1)}{2} + 2$. We will consider two cases.

\medskip\textit{Case 1.} $a_{1} \geq 0$. When $d_{0} > 0$, the $x_{0}^{d_{0}}x_{1}^{d_{1}}$ monomial has $a_{0}d_{0} + a_{1}d_{1} > a_{1}$, so it is zero for all
$q_{i}$'s except $q_{0}$; when $d_{0} > 1$ it is also zero for $q_{0}$, since $a_{0}d_{0} + a_{1}d_{1} \geq 2a_{0}$. This gives us a total codimension of $n^{2} +
(n-1)$, which is larger than $\frac{n(n+1)}{2} + 1$ for all $n \geq 2$. When $n = 1$ this case is impossible because we need to have $a_{0} + a_{1} = 0$.

\medskip\textit{Case 2.} $a_{1} < 0$. We have $a_{0} = -(a_{1} + \ldots + a_{n}) > -a_{i}$ for all $i$; therefore, the $x_{0}^{d-1}x_{i}$ monomial is zero in every
$q_{j}$ except $q_{0}$; the $x_{0}^{d}$ monomial is always zero. This gives us a codimension of $n^{2} + n + 1$, which is large enough for all $n$.\end{proof}

\begin{rem}The larger spaces $\Hom_{d}^{n, s}$ and $\Hom_{d}^{n, ss}$ have a meaning in the field of moduli spaces more than in this of dynamical systems, where we study
the iterates of morphisms. The problem is that we cannot always iterate rational maps which are not morphisms, even if they are stable: the image may not be dense,
and may eventually map to a locus on which the map is ill-defined. A map of the form $(q:0:0:\ldots:0)$ with $q(1, 0, \ldots, 0) = 0$ will be impossible to iterate.
For general $q$, it will also be stable for large $d$, because we will have $a_{0}d_{0} + \ldots + a_{n}d_{n} > a_{0}$ for many different $\mathbf{d}$'s no matter
how we choose the $a_{i}$'s, even after conjugation. When $n = 1$, it suffices to have $d \geq 4$, because then $\varphi$ is unstable only if is of the form $(p:q)$
with $p$ and $q$ sharing a common root of multiplicity at least $\frac{d-1}{2}$, and we may pick a map $(q:0)$ with $q$ having distinct roots. For one approach for
giving a completion of $\Hom_{d}^{n}$ in a way that permits iteration at the boundary, see \cite{DeM05}.\end{rem}

\section{Stabilizer Groups}

\medskip The moduli space $\mathrm{M}_{d}^{n}$, as well as its stable and semistable completions, has a well-defined function mapping each morphism to its
stabilizer group in $\PGL(n+1)$, which will be well-defined up to conjugation. This stabilizer will be finite, at least on $\mathrm{M}_{d}^{n, s}$, from standard
facts from geometric invariant theory. We will study the possible subgroups of $\PGL(n+1)$ that may occur as stabilizers of morphisms. We gain very little by
assuming Theorem~\ref{stable}, so we might as well not assume it \emph{a priori}; this will provide an alternative proof for it.

\medskip Note that the resultant is a $\PGL(n+1)$-invariant section of a $\PGL(n+1)$-linearizable divisor on $\mathbb{P}^{N}$ that is nonzero on $\Hom_{d}^{n}$.
Therefore, on $\Hom_{d}^{n}$ stability is equivalent to having closed fibers, which is equivalent to having a stabilizer group of the lowest possible dimension (see
chapter $1$ of \cite{GIT}). Hence, to provide a second proof of Theorem~\ref{stable}, it suffices to show that the stabilizer of every $\varphi \in \Hom_{d}^{n}$ is
finite. This was done in \cite{PST}. We will prove a stronger result:

\begin{thm}\label{finite}The stabilizer of every point in $\Hom_{d}^{n}, d > 1$, is a finite group of order bounded in terms of $n$ and $d$.\end{thm}

\begin{proof}Note that if $A \in \Stab(\varphi)$, then $BAB^{-1} \in \Stab(B\varphi B^{-1})$. Therefore, when considering individual stabilizing matrices, we may
assume they are in Jordan canonical form. We use the following result:

\begin{lem}\label{diagonal}If $A \in \Stab(\varphi)$, and $\varphi$ is not purely inseparable, then $A$ is diagonalizable.\end{lem}

\begin{proof} In characteristic zero, this is trivial given Theorem~\ref{stable}. However, it is not trivial in characteristic $p$; the proof works for every
characteristic, so we lose nothing from not using Theorem~\ref{stable}.

\medskip We will assume that $A$ is not diagonalizable and derive a contradiction. It suffices to assume that $A$ is a Jordan matrix whose largest Jordan block is of
size $r > 1$. After conjugation and scaling, we may assume that the first Jordan block is also the largest, and has eigenvalue $1$. We will label the rows and
columns from $0$ to $n$, in parallel with the labels for the $q_{i}$'s. We will also write $\varphi = (q_{0}:q_{1}:\ldots:q_{n})$, $k_{i} = a_{ii}$ for the
eigenvalue in the $i$th position, and $r_{i}$ for the size of the Jordan block containing $a_{ii}$. We have $r_{0} = r, k_{0} = 1, r_{i} \leq r$.

\medskip Note that the inverse of the first Jordan block is the matrix with zeroes below the main diagonal and $a_{ij} = (-1)^{i-j}$ on or above it. Therefore, each
vector $\mathbf{x} = (x_{0}, x_{1}, \ldots, x_{n})$ is transformed to:

$$\mathbf{x'} = (x_{0}-x_{1}+\ldots\pm x_{r-1}, x_{1}-x_{2}+\ldots\mp x_{r-1}, \ldots, x_{r-1}, \ldots, \frac{1}{k_{n}}x_{n})$$

\noindent We write $q_{i}'(\mathbf{x}) = q_{i}(\mathbf{x'})$. Similarly, $A$ transforms $\varphi = (q_{0}, \ldots, q_{n})$ to:

$$\varphi' = (q_{0}' + q_{1}', q_{1}' + q_{2}', \ldots, q_{r-1}', \ldots, k_{n}q_{n}')$$

\noindent Since $A$ stabilizes $\varphi$, we need $\varphi'$ to be a scalar multiple of $\varphi$.

\medskip For each $\mathbf{d} \in \mathbb{Z}^{n+1}$, we denote the $\mathbf{x^{d}}$ coefficient of $q_{i}$ (respectively $q_{i}'$) by $c_{\mathbf{d}}(i)$ (resp.
$c_{\mathbf{d}}'(i)$). We suppress trailing zeroes for simplicity, so that $c_{d}$ denotes the $x_{0}^{d}$ coefficient. We are looking for the largest $i$ such that
$c_{d}(i) \neq 0$; such an $i$ exists, or else $(1:0:\ldots:0)$ is a common root of all the $q_{i}$'s. As the only $x_{0}^{d}$ term in $\mathbf{x'^{d}}$ comes from
$x_{0}'^{d}$, we have $c_{d}'(j) = c_{d}(j)$ for all $j$. Now in $\varphi'$, the $i$th term is either $q_{i}'$ or $q_{i}' + k_{i}q_{i+1}'$, so that its $x_{0}^{d}$
coefficient is $k_{i}c_{d}(i)$. This implies that the scaling factor is $k_{i}$, i.e. $\varphi' = k_{i}\varphi$.

\medskip Now, assume that $i$ is not at the beginning of its Jordan block, that is that $a_{i-1, i} = 1$. Then $k_{i-1} = k_{i}$, and the fact that $\varphi' =
k_{i}\varphi$ implies that $k_{i-1}c_{d}'(i-1) + c_{d}'(i) = k_{i}c_{d}(i-1)$. This reduces to $c_{d}(i) = 0$, a contradiction. Therefore, $i$ is at the beginning
of its Jordan block.

\medskip Let us now consider the $x_{0}^{d-1}x_{1}$ coefficients, and assume throughout that all indices are in the same Jordan block as $i$. We have $c_{d-1,
1}'(j) = c_{d-1, 1}(j) - dc_{d}(j)$. For $j > i$, this reduces to $c_{d-1, 1}'(j) = c_{d-1, 1}(j)$. Conversely, the corresponding term to $c_{d-1, 1}$ in $\varphi'
= k_{i}\varphi$ will be $k_{i}c_{d-1, 1}'(j) + c_{d-1, 1}'(j+1) = k_{i}c_{d-1, 1}(j)$. When $j > i$, this implies that $c_{d-1, 1}'(j+1) = 0$, so that $c_{d-1,
1}(j) = 0$ for $j > i + 1$; conversely, for $i + 1$, we obtain $k_{i}c_{d-1, 1}'(i) + c_{d-1, 1}'(i+1) = k_{i}c_{d-1, 1}(i)$, which reduces to $c_{d-1, 1}(i+1) =
k_{i}dc_{d}(i) \neq 0$. This shows that $i + 1$ is the largest index with a nonzero $x_{0}^{d-1}x_{1}$ coefficient, at least in the Jordan block containing $i$.

\medskip We may apply induction on $s(\mathbf{d}) = d_{1} + 2d_{2} + \ldots + (r-1)d_{r-1}$, and find that in the Jordan block containing $i$, the largest index
with a nonzero $\mathbf{x^{d}}$ coefficient is $i + s(\mathbf{d})$. Note that the Jordan block has $r_{i} \leq r$ elements, but the number of monomial indices
attached to the first Jordan block is $(r-1)d + 1$, which is strictly greater than $r$ when $d, r > 1$. This is a contradiction: the last element of the Jordan
block has $k_{i}c_{\mathbf{d}}' = k_{i}c_{\mathbf{d}}$ for all $\mathbf{d}$, i.e. $c_{\mathbf{d}}(i + r_{i} - 1)' = c_{\mathbf{d}}(i + r_{i} - 1)$, but that last
equality is only true when $s(\mathbf{d}) \leq r_{i}$, which is not the case for all $\mathbf{d}$. Since we are assuming $d > 1$, we must have $r = 1$, and we are
done.

\medskip The careful reader may note that the proof that $i + s$ is the largest index with a nonzero $\mathbf{x^{d}}$ coefficient for $s(\mathbf{d}) = s$ makes an
assumption about the characteristic we are working in. In characteristic zero, $d \neq 0$ and there is no problem. In characteristic $p$, we need to treat
separately the case when $p < d$. Then for example we may have $p \mid d$, so that $c_{d-1, 1}'(j) = c_{d-1, 1}(j)$ for all $j$, and $c_{d-1, 1}(i+1)$ may be zero.
Note that the number of monomial indices containing $x_{0}^{d-2}$ attached to the first Jordan block is $2(r-1) + 1$, which is strictly greater than $r$ when $r >
1$; when $p \nmid d(d-1)$, we may restrict ourselves to such monomials, and the proof proceeds as in characteristic zero.

\medskip When $p \mid d-1$, we may restrict ourselves to monomials containing $x_{0}^{d-1}$, and proceed with the proof. We will only encounter an
obstruction if $r_{i} = r$ and only at the end of the Jordan block, where the existence of a nonzero $x_{0}^{d-1}x_{r-1}$ monomial does not guarantee that of
$x_{0}^{d-2}x_{1}x_{r-1}$. However, the action of $A$ on $q_{i + r - 1}$ takes it to $k_{i}q'_{i + r - 1}$, and we must have $c_{\mathbf{d}}'(i + r - 1) =
c_{\mathbf{d}}(i + r - 1)$ for all $\mathbf{d}$. If we write $d-1 = p^{l}m, m \nmid p$, then we see that $x_{0}^{d-1}x_{r-1}$ is transformed to $k_{i}(x_{0} - x_{1}
+ \ldots \pm x_{r-1})^{d-1}x_{r-1} = k_{i}(x_{0}^{p^{l}} - \ldots \pm x_{r-1}^{p^{l}})^{m}x_{r-1}$ which shows that the $x_{0}^{p^{l}(m-1)}x_{1}^{p^{l}}$ monomial
does not satisfy $c_{\mathbf{d}}'(i + r - 1) = c_{\mathbf{d}}(i + r - 1)$. This yields a contradiction.

\medskip Finally, when $p \mid d$, we may write $d = p^{l}m$. When $m > 1$, we apply exactly the same proof as in characteristic zero, except that we write $m$
instead of $d$ and $m_{j} = \frac{d_{j}}{p^{l}}$ instead of $d_{j}$; then we define $s(\mathbf{d}) = m_{1} + \ldots + (r-1)m_{r-1}$, and in the Jordan block
containing $i$, the largest index with a nonzero $\mathbf{x^{d}}$ coefficient is $i + s(\mathbf{d})$. As $m > 1$, we have $(r-1)m + 1 > r$ for $r > 1$, and we have
the same contradiction as in the characteristic zero case. Note that when $m = 1$, we may derive the same contradiction from any nonzero monomial not of the form
$x_{j}^{d}$, which must exist if $\varphi$ is not purely inseparable. Hence, if $\varphi$ has a non-diagonalizable stabilizer then it is purely inseparable and we
are done.\end{proof}

\medskip With the above lemma, we know that any abelian subgroup of $\Stab(\varphi) \in \GL(n+1)$ will be simultaneously diagonalizable. We will prove the following
uniform bound on the size of abelian stablizing subgroups:

\begin{lem}\label{bound}Every diagonal subgroup stabilizing $\varphi \in \Hom_{d}^{n}$ is of size at most $d^{n+1}$.\end{lem}

\begin{proof}A diagonal matrix $A$ with diagonal entries $(a_{0}, a_{1}, \ldots, a_{n})$ acts on each $q_{i}$ by multiplying $c_{\mathbf{d}}(i)$ by
$\frac{a_{i}}{a_{0}^{d_{0}}\ldots a_{n}^{d_{n}}}$. Our case of interest will be the $x_{j}^{d}$ coefficients. Each has to be nonzero for at least one $i$, which
induces the equation $a_{i}$ = $a_{j}^{d}$. Note that we may set the scaling factor $k$ to be $1$, since the scalar matrix $k^{\frac{1}{1-d}}$ multiplies every
coefficient by $k$.

\medskip Now, we have at least $n + 1$ different relations $a_{i} = a_{j}^{d}$. We may drop relations until each $j$ has just one $i$ such that such a relation
holds; dropping relations will increase the size of the group, so by bounding the size of the larger group, we will bound the size of any automorphism group.

\medskip We obtain a function $j \mapsto i$. If the function is bijective, we may write it as a product of disjoint cycles, and conjugate to get the cycles to be
$(0\ 1\ \ldots\ s_{1}-1) \ldots (n-s_{k}+1\ \ldots\ n)$, where here $r_{i}$ denotes the length of the $i$th cycle, and has nothing to do with the definition in
Lemma~\ref{diagonal}. Then $a_{0}^{d^{r_{1}}} = a_{0}$ and $a_{0}$ is a root of unity of order dividing $d^{r_{1}} - 1$, the choice of which uniquely determines
$a_{i}, 0 \leq i \leq r_{1} - 1$. We have similar results for $a_{r_{1}}, \ldots, a_{n-r_{k}+1}$; since $\sum r_{i} = n+1$, this bounds the size of the group by
$d^{n+1}$.

\medskip In general, of course, the function $j \mapsto i$ may not be bijective, so we can only write it as a product of precycles, whose cycles are disjoint. Here
a precycle means a cycle and zero or more tails. The above discussion applies to the cycles. For the tails, suppose without loss of generality that $(0\ 1\ \ldots\
r)$ is a tail where $r$ and no element before it is part of a cycle; then the choice of $a_{r}$ determines a choice of $d$ possibilities for $a_{r-1}$ and in
general $d^{s}$ for $a_{r-s}$ subject to the obvious compatibility condition. This clearly respects the bound of $d^{n+1}$: if $m$ is the total number of elements
in cycles, then we have at most $m^{n+1}$ possibilities for the cycles, each of which gives us exactly $(d-m)^{n+1}$ possibilities for the tails.\end{proof}

\medskip The bound $d^{n+1}$ works for abelian stabilizing subgroups in the purely inseparable case as well. We may view a purely inseparable $\varphi$ as the action
of raising every coefficient to the $d$th power followed by the matrix $B$. Then $A\varphi A^{-1} = \varphi$ if and only if $ABA^{-1}_{d} = B$, where $A_{d}$ is the
image of the matrix $A$ under the homomorphism of raising every entry to the $d$th power; we need to show the group of such $A$, which we will write as $\Stab(B)$,
is finite. Since $A$ and $A_{d}$ are conjugate, all eigenvalues of $A$ are in $\mathbb{F}_{d}$.

\medskip We may conjugate an abelian stabilizing subgroup $G$ to obtain a block diagonal group with each block upper triangular and with its $(i, j)$ entry depending
only on $j - i$. We may also fix one element, $C$ to be in Jordan canonical form, in which case we will have $C_{d} = C$ and thus $BC = CB$. Then $B$ is in block
form; labeling the blocks by $r, s$ and the $r$th block of $C$ by $C_{r}$, we see that the $B_{rs}$ is nonzero if and only if the blocks $r$ and $s$ are of the same
size and equal for every element of $\Stab(B)$, and in any case $B_{rs}$ commutes with $C_{r} = C_{s}$, so it is upper triangular with its $(i, j)$ entry depending
only on $j - i$. In particular, it commutes with every $A_{r} = A_{s}$, so that $B$ commutes with $G$. Hence for all $A \in G$, we have $AB = BA$ and $ABA^{-1}_{d}
= B$, so that $A = A^{d}$ and $A$ has entries in $\mathbb{F}_{d}$. Furthermore, for each block in $G$ of size $r$, we have $r$ positive possibilities for $j - i$,
inducing $d^{r}$ possible blocks, and $d^{n+1}$ possible matrices in $G$.

\medskip Note that we may have additional stabilizing matrices in $\PGL(n+1)$. These occur when there exists an automorphism of the set $\{0, 1, \ldots, n\}$ that
does not leave the diagonal vector $\mathbf{a} = (a_{0}, \ldots, a_{n}) \in \mathbb{A}^{n+1}$ fixed, but does fix $\mathbf{a} = (a_{0}:\ldots:a_{n}) \in
\mathbb{P}^{n}$. Since the automorphism has to fix $a_{0}a_{1}\ldots a_{n}$, we see that it must send each $a_{i}$ to $\zeta a_{i}$ where $\zeta$ is a root of unity
of order at most $n+1$; hence there are at most $n+1$ possibilities for such an automorphism, modulo automorphisms that fix $\mathbf{a} \in \mathbb{A}^{n+1}$ and
are hence simultaneously block-diagonalizable with $A$.

\medskip We will rely on one final bound, due to G. A. Miller \cite{Mil}:

\begin{prop}\label{group}The size of a finite group is bounded in terms of the size of its largest abelian subgroup.\end{prop}

\begin{proof}It suffices to show this for $p$-groups. For each $n$, we let $k(n)$ be the minimal exponent of the largest abelian subgroup of any $p$-group of
exponent $n$. Furthermore, for each $l \leq n$, we let $k(n, l)$ be the minimal exponent subject to the restriction that $Z = Z(G)$ have exponent $l$, so that $k(n)
= \min\{k(n, l)\}$. It is enough to show that $\lim_{n \to \infty}k(n) = \infty$.

\medskip It is trivial to show that $k(2) = 2$. In general, for a $p$-group of exponent $n$ and center of exponent $l$, let $g$ be such that $g \notin Z$,
$g^{p} \in Z$, and $gZ \in Z(G/Z)$. Unless $G$ is abelian, in which case the result is trivial, we may take $g$ to be a preimage of a nontrivial element in the
socle of $G/Z$. For every $h \in G$, $hgh^{-1} = gz$ for some $z \in Z$; we obtain a group homomorphism $h \mapsto z$ from $G$ to $Z$. The homomorphism has kernel
$K$ of exponent at least $n-l$ and center containing $\langle Z, g\rangle$. Any abelian subgroup of $K$ will be an abelian subgroup of $G$, so that we obtain $k(n,
l) \geq k(n-l, l+1)$. It easily follows that $k(n) \geq 2\sqrt{n}$.\end{proof}

\medskip The bound in the above proposition is very weak. It is known that for odd $p$ we have $k(n) \leq \frac{n+4}{3}$ and for $p = 2$ we have $k(n) \leq
2\frac{n+3}{5}$ \cite{Dan}, but little more. It is also not known \emph{a priori} that the group has to be finite, only that if it is finite then it is bounded. We
may use Theorem~\ref{stable} and finish. However, with little additional effort, we may prove finiteness directly, providing an alternative proof that all morphisms
are stable. The fact that finite implies uniformly bounded means that it is enough to show that every finitely generated stabilizing subgroup is finite. More
precisely:

\begin{prop}\label{inject}Every finitely generated subgroup of $\PGL(n)$ contained in finitely many finite-order conjugacy classes is finite.\end{prop}

\begin{proof}Let $R$ be the $\mathbb{Z}$-algebra generated by the finitely many coefficients of the generators. Then the group is contained in $\PGL(n, R)$, and we
may project it into the finite group $\PGL(n, R/\mathfrak{m})$ where $\mathfrak{m}$ is a maximal ideal in $R$; we will show the map can be chosen to be injective.
In fact, each non-unipotent conjugacy class $i$ contains two different eigenvalues, $a_{i_{1}}, a_{i_{2}}$; therefore, if we choose $\mathfrak{m}$ not to contain
$a_{i_{1}} - a_{i_{2}}$, which we can since there are only finitely many such elements, then the map will have unipotent kernel. In characteristic $0$, the only
finite-order unipotent matrix is the identity, so the map is injective and we are done.

\medskip In characteristic $p$, we obtain a finite-index and hence finitely generated unipotent group. We may conjugate it by some matrix $P$ to be upper
triangular; then matrix multiplication is equivalent to addition of the $(r, r+1)$ entry for any $r$, and the finite generation implies that the set of all $(r,
r+1)$ entries lies in a finitely generated $\mathbb{Z}/p\mathbb{Z}$-vector space, which is finite. For the matrices with all $(r, r+k)$ entries for all $k \leq l$,
matrix multiplication corresponds to addition of $(r, r+l+1)$ entries, and we may add those entries to our vector space, which will remain finite. We may now
construct $\mathfrak{m}$ to avoid the finite vector space and the determinant of $P$, as well as the eigenvalue differences described above. The map will then be
injective.\end{proof}

\medskip Note that in the proof of proposition we make no assumption on the base ring. Of course, the argument in the proposition applies to $\GL(n+1)$, and shows
that the answer to Burnside's problem, which asks whether a finitely generated group of bounded exponent is necessarily finite, is yes when restricted to subgroups
with faithful finite-dimensional representations over any field.\end{proof}

\medskip For each stabilizer group $G \in \PGL(n+1)$, there is a closed subscheme $\Fix(G) \in \Hom_{d}^{n}$ consisting of all $\varphi$ with stabilizer group
containing $G$. Theorem~\ref{finite} states that every $G$ with nonempty $\Fix(G)$ is finite and of bounded order. Furthermore, each nontrivial stabilizing matrix
is, up to conjugation, one of the $d^{n+1}$ possibilities for each of the $(n+1)^{n+1}$ functions on the set $\{0, 1, \ldots, n\}$. We may strengthen this result as
follows:

\begin{cor}\label{stabilizer}There are only finitely many $G$ with nonempty $\Fix(G)$ up to conjugation. In particular, on an open dense set of $\Hom_{d}^{n}$,
which descends to $\mathrm{M}_{d}^{n}$, the stabilizer group is trivial.\end{cor}

\begin{rem}The statement that there are only finitely many such $G$ up to conjugation is stronger than the statement that there are only finitely many $G$ up to
isomorphism, which follows trivially from the bound on the size of $G$.\end{rem}

\begin{proof}Since the size of $G$ is bounded, it suffices to show that each stabilizing subgroup has finitely many projective $n+1$-dimensional representations up
to conjugacy. This is always true when the representation is completely reducible, which will be true if the ambient characteristic $p$ does not divide
$\left|G\right|$. But when $\Fix(G)$ is not purely inseparable, every element will be diagonalizable, so it will have order not divisible by $p$, so that $G$ has
order not divisible by $p$. In the purely inseparable case, we have $\PGL(n+1)$ acting on itself stably and with finite stabilizers, so that each orbit is of
dimension $(n+1)^{2} - 1$ and thus consists of all of $\PGL(n+1)$. In other words, every purely inseparable map is, up to conjugation,
$(x_{0}^{d}:\ldots:x_{n}^{d})$, so that its stabilizer group is conjugate to $\PGL(n+1, \mathbb{F}_{d})$.

\medskip It remains to be shown that the complement of $\bigcup_{G \supset I}\Fix(G)$ is dense; its openness follows from the fact that the
condition $A\varphi A^{-1} = \varphi$ is closed. It suffices to show that each $\Fix(G)$ is a proper subset of $\Hom_{d}^{n}$. We lose nothing if we ignore purely
separable maps. From the proof of Lemma~\ref{bound}, each of the finitely many elements that may occur in $G$, a diagonal matrix with $i$th entry $a_{i}$,
multiplies $c_{\mathbf{d}}(i)$ by $\frac{a_{i}}{\mathbf{a^{d}}}$, and hence induces the relation $c_{\mathbf{d}}(i) = 0$ outside a set of $(\mathbf{d}, i)$'s for
which $\frac{a_{i}}{\mathbf{a^{d}}}$ is constant. If $\frac{a_{i}}{\mathbf{a^{d}}}$ is constant for all $(\mathbf{d}, i)$, then we have $a_{i} = k\mathbf{a^{d}}$;
choosing a constant $\mathbf{d}$, we see that $a_{i}$ is constant, so $A$ is a scalar matrix. Hence no non-trivial $A$ fixes all of $\Hom_{d}^{n}$.\end{proof}

\medskip Note that when $n = 1$, \cite{Sil96} has an explicit bound on the size of $\Stab(\varphi)$ of $n_{1}!n_{2}!n_{3}!$, where the $n_{i}$'s are indices for
which there exist periodic points for $\varphi$ of exact order $n_{i}$. The technique in this paper improves on that bound. Following the proof of
Lemma~\ref{bound}, we have three possibilities for the map $j \mapsto i$ up to conjugation: $(1, 2) \mapsto (1, 2)$, $(1, 2) \mapsto (2, 1)$, and $(1, 2) \mapsto
(1, 1)$. In the first case, $a_{0} = \zeta_{d-1}^{i}$ and $a_{1} = \zeta_{d-1}^{j}$, where we use $\zeta_{i}$ to denote an $i$th root of unity; modulo multiplying
both $a_{0}$ and $a_{1}$ by some $\zeta_{d-1}$, we obtain a cyclic group of order $d-1$. In the second case, we have $a_{0} = \zeta_{d^{2} - 1}^{i}, a_{1} =
a_{0}^{d}$, and modulo multiplying both by $\zeta_{d^{2} - 1}^{d+1}$, we obtain a cyclic group of order $d+1$. In the third case, $a_{0} = \zeta_{d-1}$ and
$a_{1}^{d} = a_{0}$, and modulo multiplying both by $\zeta_{d-1}$, we obtain a cyclic group of order $d$.

\medskip Thus every diagonalizable abelian subgroup $A$ of $\Stab(\varphi)$ will be cyclic of size dividing $d-1$, $d$, or $d+1$. Furthermore, the only
non-diagonalizable element commuting with $A$ can be the matrix $M$ corresponding to the automorphism permuting $x_{0}$ and $x_{1}$; we have $M^{-1} = M$ and $MAM =
A$ in $\PGL(2)$ if and only if $\frac{a_{1}}{a_{0}} = \frac{a_{0}}{a_{1}}$, or, equivalently, $a_{i} = \pm 1$ for $i = 0, 1$. In other words, the only possible
non-diagonalizable abelian subgroup $A$ is $\mathbb{Z}/2\mathbb{Z} \times \mathbb{Z}/2\mathbb{Z}$.

\medskip Now, the only finite subgroups of $\PGL(2)$ are, up to conjugation, cyclic, dihedral, tetrahedral, octahedral, or icosahedral \cite{Sil94}. The last three
groups are of order at most $60$; only the first two are infinite families. Since the largest abelian subgroup of the dihedral group of order $2k$ is of order $k$,
we see that for large $d$, the order of $\Stab(\varphi)$ is bounded by $2(d+1)$.

\medskip We conclude this section with a remark that $\mathrm{M}_{d}^{n}(k)$, consisting of all $k$-rational points in $\mathrm{M}_{d}^{n}(\overline{k})$, is not
the same as the quotient $\Hom_{d}^{n}(k)/\PGL(n+1, k)$. The latter parametrizes morphisms of $\mathbb{P}^{n}_{k}$ up to conjugation defined over $k$, the former up
to conjugation defined over $\overline{k}$. There exist maps defined over $k$ which are conjugate over $\overline{k}$ but not over $k$ itself. For examples, see
\cite{Sil96} and \S\S 4.7-4.10 of \cite{ADS}.

\section{Rationality of $\mathrm{M}_{d}$}

\noindent In this section, we show that when $n = 1$, the variety $\mathrm{M}_{d} = \mathrm{M}_{d}^{1}$ is rational. This partly generalizes Silverman's result in
\cite{Sil96} that $\mathrm{M}_{2} = \mathbb{A}^{2}$ over $\mathbb{Z}$. We do so by parametrizing fixed points of $\varphi$. The fixed point set of $\varphi$,
$\Fix(\varphi)$, is the intersection of two curves in $\mathbb{P}^{1} \times \mathbb{P}^{1}$, the graph $\Gamma_{\varphi}$ and the diagonal embedding $\Delta$. As
$\Delta$ is irreducible and not contained in $\Gamma_{\varphi}$ for $d > 1$, this is a proper intersection of divisors of type $(1, 1)$ and $(d, 1)$, so it has $d +
1$ points, counting multiplicity. We have:

\begin{thm}\label{rational}$\mathrm{M}_{d}$ is birational to the total space of a rank-$d$ vector bundle on $\mathrm{M}_{0, d+1}$, the space of unmarked $d+1$
points on $\mathbb{P}^{1}$. Since $\mathrm{M}_{0, d+1}$ is rational, it follows that $\mathrm{M}_{d}$ is rational.\end{thm}

\begin{proof}We explicitly write $\varphi(x:y) = (p:q)$ where $p(x, y) = a_{d}x^{d} + \ldots + a_{0}y^{d}$ and $q(x, y) = b_{d}x^{d} + \ldots + b_{0}y^{d}$. The
fixed points of $\varphi$ are those for which $(p:q) = (x:y)$, which are the roots of the homogeneous degree-$d + 1$ polynomial $py - qx$. The polynomial $py - qx$
induces a map from $\Rat_{d}$ to $(\mathbb{P}^{1})^{d+1}/S_{d+1}$ where $S_{d+1}$ acts by permutation of the factors. We will call this map $\Fix$. We use the
following lemma:

\begin{lem}\label{surjective}The map $\Fix$ is surjective, and has rational fibers.\end{lem}

\begin{proof}A point $(x:y)$ is fixed if and only if we have $py = qx$, i.e. $a_{d}x^{d}y + \ldots + a_{0}y^{d+1} = b_{d}x^{d+1} + \ldots + b_{0}xy^{d}$. This is a
homogeneous linear condition in the coefficients of $\varphi$, and we have $d+1$ such conditions compared with $2d + 2$ variables. From elementary linear algebra,
we have a solution space of linear dimension $d+1$, or projective dimension $d$. It is a linear subvariety of $\mathbb{P}^{2d+1}$, so it is rational.

\medskip We can also show that this dimension-$d$ space will not be contained in the resultant locus. We fix a set of fixed points and write $r$ for the polynomial
having those fixed points as roots. We need to show $r$ is of the form $py - qx$ for some $p$ and $q$ sharing no common root. By conjugating, we may assume neither
$(0:1)$ nor $(1:0)$ is a root of $r$, so that it has a nonzero $x^{d+1}$ coefficient, which we may take to be $1$, and a nonzero $y^{d+1}$ coefficient. Now we let
$q = -x^{d}$ so that $r + qx$ is divisible by $y$, yielding $p = \frac{r + qx}{y}$. Now $r + qx$ has a nonzero $y^{d+1}$ coefficient, so $p$ has a nonzero $y^{d}$
coefficient; therefore, $p$ does not have $(0:1)$ as a root, so it shares no root with $y$.\end{proof}

\medskip Now, $\Fix$ descends to a rational map $\Fix':\mathrm{M}_{d} \to (\mathbb{P}^{1})^{d+1}/S_{d+1}\PGL(2)$ where $\PGL(2)$ acts diagonally; we are restricting
to the open set of $\mathrm{M}_{d}$ whose fixed points are in the stable space of the action of $\PGL(2)$ on $(\mathbb{P}^{1})^{d+1}/S_{d+1}$. With this
restriction, the image is $\mathrm{M}_{0, d+1}$, so it suffices to show the general fiber of $\Fix'$ is rational. Lemma~\ref{surjective} says that the fiber of
$\Fix$ is rational, so it suffices to show that the automorphism group of the general point in $(\mathbb{P}^{1})^{d+1}/S_{d+1}$ is small enough that the quotient of
the fiber by it is still rational. Using Noether's problem \cite{Noe} \cite{Sal}, we will show a stabilizer of size $4$ or $6$ is small enough.

\begin{lem}\label{auto}Let $d > 1$. The automorphism group of a general configuration of $d + 1$ unmarked points in $\mathbb{P}^{1}$ is trivial, unless $d = 2$, in
which case it is $S_{3}$, or $d = 3$, in which case it is $\mathbb{Z}/2\mathbb{Z} \times \mathbb{Z}/2\mathbb{Z}$.\end{lem}

\begin{proof}We will use inhomogeneous coordinates. For $d = 2$, we can conjugate the three points to be $0, 1, \infty$; the set is then stabilized by every
permutation in $S_{3}$, so it has size $6$. For $d > 3$, we will show that the stabilizer is generically trivial, and on the way show that for $d = 3$ the
stabilizer is generically of order $4$, consisting of all elements in $S_{4}$ of cycle type $(2, 2)$. This will be enough to prove the theorem.

\medskip First, note that if a $(d+1)$-cycle stabilizes the set of points, then by conjugation we may assume it sends $0$ to $1$, $1$ to $\lambda$, $\mu$ to
$\infty$, and $\infty$ to $0$. The cycle, regarded as an element of $\PGL(2)$, is of the form $\frac{ax + b}{cx + e}$; then $\frac{b}{e} = 1$, $a = 0$, $\frac{a +
b}{c + e} = \lambda$, and $c\mu + e = 0$. These equations together imply that $\lambda = \frac{e}{c + e} = \frac{e}{e - \frac{e}{\mu}} = \frac{\mu}{\mu - 1}$. For a
generic choice of $\mu, \lambda$, this can never happen, so no $(d+1)$-cycle is in the stabilizer. This remains true for $d = 3$, in which case we are forced to
have $\lambda = \mu$, since generically $\lambda \neq \frac{\lambda}{\lambda - 1}$.

\medskip Observe that if an automorphism of cycle type $(c_{1}, \ldots, c_{k})$ stabilizes the set, then each subset corresponding to the $i$th cycle is stabilized
by a $c_{i}$-cycle. Therefore, the above discussion shows that no cycle of length $4$ or more stabilizes a generic set. We have reduced to the case when all cycles
are of size $1$, $2$, or $3$. Now, if we have a stabilizing automorphism which includes a $3$-cycle, we may conjugate the $3$-cycle to be $(0\ 1\ \infty)$, forcing
it to act on $\mathbb{P}^{1}$ as $\frac{1}{1-x}$. Generically, if $\lambda$ is a fourth point, none of the points in the set (including $\lambda$) will be
$\frac{1}{1 - \lambda}$. We are left with cycles of size $1$ or $2$. If we have a stabilizing automorphism with two $2$-cycles, then up to conjugation we may assume
the element acts on four points as $(0\ \infty)(1\ \lambda)$, so that it maps $x$ to $\frac{\lambda}{x}$. If $d = 3$ then this will stabilize the set regardless of
what $\lambda$ is. If $d > 3$ then we have an additional point $\mu$, and generically $\frac{\lambda}{\mu}$ will not be in our set.

\medskip We are left with automorphisms that act as single $2$-cycles, fixing $d - 1$ points. For $d \geq 4$, they will fix $3$ points and therefore act trivially.
For $d = 3$, we may assume by conjugation that the element acts as $(0\ 1)$ and fixes $\infty$; this forces it to be the automorphism $1 - x$, which generically
does not fix $\lambda$. This leaves us with automorphisms consisting only of $1$-cycles, i.e. the identity.\end{proof}

\medskip We will return to Noether's problem now. Let us work over a fixed field $k$. Recall \cite{Noe} that if $K = k(x_{1}, \ldots, x_{m})$ is a purely
transcendental field, and $G$ is a finite group of size $2$, $3$, $4$, or $6$ permuting the $x_{i}$'s, then $K^{G}$ is purely transcendental as well. In particular,
if $R$ is the graded $k$-algebra $k[x_{1}, \ldots, x_{m}]$, and $G$ acts on it by permutation of the $x_{i}$'s, then $\Proj R^{G}$ is rational. We will show this to
be the case when $R$ is the fiber of $\Fix$ in the $d = 2$ and $d = 3$ cases, by finding an orbit $y_{1}, \ldots, y_{m}$ generating $R$ over $k$.

\medskip When $d = 2$, we have a $2$-dimensional fiber. Explicitly, we have six homogeneous variables $a_{i}, b_{i}, 0 \leq i \leq 2$, on which the automorphism
group $\PGL(2)$ acts linearly. The fiber we are interested in consists of maps fixing the points $0, 1, \infty$, corresponding to the linear conditions $a_{0} = 0$,
$a_{0} + a_{1} + a_{2} = b_{0} + b_{1} + b_{2}$, $b_{2} = 0$, respectively. The values of $a_{2}, a_{1}, b_{0}$ uniquely determine that of $b_{1}$, so we may write
the fiber as $\Proj k[a_{2}, a_{1}, b_{0}]$. The group $S_{3}$ acts linearly and faithfully on the $k$-vector space spanned by $a_{2}, a_{1}, b_{0}$. Let us
consider the action of the automorphism $(0\ \infty) = \frac{1}{x}$:

$$\varphi(x) = \frac{a_{2}x^{2} + a_{1}x}{b_{1}x + b_{0}}$$
$$\frac{1}{\varphi(\frac{1}{x})} = \frac{b_{0}x^{2} + b_{1}x}{a_{1}x + a_{2}}$$
$$a_{2} \mapsto b_{0}$$
$$a_{1} \mapsto b_{1} = a_{2} + a_{1} - b_{0}$$
$$b_{0} \mapsto a_{2}$$

\noindent Observe that this automorphism fixes $a_{2} + b_{0}$. Let us also consider the action of the automorphism $(0\ 1) = 1 - x$:

\begin{align}\nonumber 1 - \varphi(1-x) & = 1 - \frac{a_{2}(1-x)^{2} + a_{1}(1-x)}{b_{1}(1-x) + b_{0}} \\
\nonumber &  = \frac{-a_{2}(1-x)^{2} + (b_{1} - a_{1})(1-x) + b_{0}}{b_{1}(1-x) + b_{0}}\end{align}
$$a_{2} \mapsto -a_{2}$$
$$a_{1} \mapsto 2a_{2} + a_{1} - b_{1} = a_{2} + b_{0}$$
$$b_{0} \mapsto b_{0} + b_{1} = a_{2} + a_{1}$$

\noindent This automorphism does not stabilize $a_{2} + b_{0}$; hence, $a_{2} + b_{0}$ has stabilizer of order $2$, and orbit of size $3$. By repeating the maps
$1-x$ and $\frac{1}{x}$, we can compute the orbit as $\{a_{2} + b_{0}, a_{1}, a_{2} + a_{1} - b_{0}\}$. This generates $R$ as long as $\chr k \neq 2$. When $\chr k
= 2$, the automorphism $1 - x$ fixes $a_{2}$, whose orbit is then $\{a_{2}, b_{0}, a_{2} + a_{1}\}$. In either case, we can construct the action of $S_{3}$ as an
action of generators, reducing the quotient to Noether's problem.

\medskip When $d = 3$, we similarly obtain a $3$-dimensional fiber, fixing the points $0, 1, \lambda, \infty$. We obtain the linear conditions $a_{0} = 0$, $b_{3} =
0$, $a_{3} + a_{2} + a_{1} = b_{2} + b_{1} + b_{0}$, $\lambda^{2}a_{3} + \lambda a_{2} + a_{1} = \lambda^{2}b_{2} + \lambda b_{1} + b_{0}$, and we may write $R$ as
$k[a_{3}, a_{2}, b_{1}, b_{0}]$. We look at the automorphism $(0\ \infty)(1\ \lambda) = \frac{\lambda}{x}$:

$$\varphi(x) = \frac{a_{3}x^{3} + a_{2}x^{2} + a_{1}x}{b_{2}x^{2} + b_{1}x + b_{0}}$$
$$\frac{\lambda}{\varphi(\frac{\lambda}{x})} = \frac{\lambda}{\frac{a_{3}\lambda^{3} + a_{2}x\lambda^{2} + a_{1}x^{2}\lambda}{b_{2}x\lambda^{2} + b_{1}x^{2}\lambda +
b_{0}x^{3}}} = \frac{b_{0}x^{3} + b_{1}\lambda x^{2} + b_{2}\lambda^{2}x}{a_{1}x^{2} + a_{2}\lambda x + a_{3}\lambda^{2}}$$
$$a_{3} \mapsto b_{0}$$
$$a_{2} \mapsto \lambda b_{1}$$
$$b_{1} \mapsto \lambda a_{2}$$
$$b_{0} \mapsto \lambda^{2} a_{3}$$

\noindent We may scale down by a factor of $\lambda$ to obtain $(\lambda^{-1}b_{0}, b_{1}, a_{2}, \lambda a_{3})$, which is equivalent to picking the representative
function $\frac{\sqrt{\lambda}}{\sqrt{\lambda^{-1}}x}$. Let us also consider the action of the automorphism $(0\ \lambda)(1\ \infty) = \frac{x - \lambda}{x - 1}$:

$$\varphi(\frac{x - \lambda}{x-1}) = \frac{a_{3}(x-\lambda)^{3} + a_{2}(x-\lambda)^{2}(x-1) + a_{1}(x-\lambda)(x-1)^{2}}{b_{2}(x-\lambda)^{2}(x-1) +
b_{1}(x-\lambda)(x-1)^{2} + b_{0}(x-1)^{3}}$$

\noindent We obtain:
$$\frac{a_{3}(x-\lambda)^{3} + (a_{2}-\lambda b_{2})(x-\lambda)^{2}(x-1) + (a_{1}-\lambda b_{1})(x-\lambda)(x-1)^{2} - \lambda b_{0}(x-1)^{3}}{a_{3}(x-\lambda)^{3}
+ (a_{2}-b_{2})(x-\lambda)^{2}(x-1) + (a_{1}-b_{1})(x-\lambda)(x-1)^{2} - b_{0}(x-1)^{3}}$$
$$a_{3} \mapsto a_{3} + a_{2} + a_{1} - \lambda(b_{2} + b_{1} + b_{0})$$

\noindent We will show the orbit of $a_{3}$ generates $R$. But first, note that $a_{3} + a_{2} + a_{1} = b_{2} + b_{1} + b_{0}$ implies that $a_{1} = b_{2} + b_{1}
+ b_{0} - a_{2} - a_{3}$, and then $\lambda^{2}a_{3} + \lambda a_{2} + a_{1} = \lambda^{2}b_{2} + \lambda b_{1} + b_{0}$ implies that $(\lambda^{2} - 1)a_{3} +
(\lambda - 1)a_{2} = (\lambda^{2} - 1)b_{2} + (\lambda - 1)b_{1}$, that is, $b_{2} = a_{3} + \frac{a_{2} - b_{1}}{\lambda + 1}$.

\medskip We have $\frac{x - \lambda}{x - 1}$ mapping $a_{3}$ to $a_{3} + a_{2} + a_{1} - \lambda(b_{2} + b_{1} + b_{0}) = (1 - \lambda)(b_{2} + b_{1} + b_{0}) = (1
- \lambda)(a_{3} + b_{0} + \frac{a_{2} + \lambda b_{1}}{\lambda + 1})$. If we then apply the map $\frac{\lambda}{x}$, we obtain $(1 - \lambda)(\lambda^{-1}b_{0} +
\lambda b_{3} + \frac{b_{1} + \lambda a_{2}}{\lambda + 1})$. The orbit is, up to scaling, $\{a_{3}, b_{0}, a_{3} + b_{0} + \frac{a_{2} + \lambda b_{1}}{\lambda +
1}, \lambda^{-1}b_{0} + \lambda b_{3} + \frac{b_{1} + \lambda a_{2}}{\lambda + 1}\}$, which generates $R$. Again, we apply Noether's problem and obtain a rational
quotient, as desired.\end{proof}

\medskip Unfortunately, this proof does not seem to generalize to $\mathrm{M}_{d}^{n}$. Although Lemma~\ref{auto} is true for all $n, d > 1$, there are two
significant obstructions. First, the dimension of the target space of the map $\Fix$ will be $n(1 + d + \ldots + d^{n})$, which is larger than $N_{d}^{n}$ unless
$n$ and $d$ are very small. This means that the map will not be surjective, though the fibers are still rational whenever they are nonempty. And second, even for
small $n$ and $d$ the base space for the vector bundle is not $\mathrm{M}_{0, d+1}$, which is relatively tame, but rather the space of $1 + d + \ldots + d^{n}$
points on $\mathbb{P}^{n}$, a much more complex object. All we can say at this stage is that $\mathrm{M}_{d}^{n}$ is unirational, which follows trivially from the
fact that it is covered by $\Hom_{d}^{n}$.

\bibliographystyle{amsplain}
\bibliography{morphisms}

\bigskip\noindent \sc{Alon Levy, Department of Mathematics, Columbia University, New York, NY 10027, USA}

\noindent \tt{email: levy@math.columbia.edu}

\end{document}